\documentclass[a4paper,10pt]{article}

\usepackage{amssymb,amsmath,graphics,mathrsfs}

\newenvironment{proofof}[1]{\par\noindent{\bf Proof of  #1.}}{\par\rightline{$\blacksquare$}}

\usepackage{amsthm}
\numberwithin{equation}{section}
\theoremstyle{definition}
\newtheorem{defi}{Definition}[section]
\newtheorem{ex}[defi]{Example}
\newtheorem*{asum*}{Assumption}

\newtheorem{rem}[defi]{Remark}
\theoremstyle{plain}
\newtheorem{lem}[defi]{Lemma}
\newtheorem{col}[defi]{Corollary}
\newtheorem{theorem}[defi]{Theorem}
\newtheorem{prop}[defi]{Proposition}

\newtheorem*{theorem*}{Theorem}
\def\be#1\ee{\begin{equation}#1\end{equation}}
\newcommand{\ba}{\begin{eqnarray} }
\newcommand{\ea}{\end{eqnarray} }
\def\bt#1\et{\begin{theorem}#1\end{theorem}}
\def\bl#1\el{\begin{lem}#1\end{lem}}
\def\bp#1\ep{\begin{prop}#1\end{prop}}
\def\bd#1\ed{\begin{def}#1\end{def}}
\DeclareMathOperator{\sgn}{sgn}

\def\ccE{{\cal E}}

\def\va{\varepsilon}
\def\ra{\rightarrow}

\def\E{\mathbf{E}}

\def\P{\mathbf{P}}

\def\R{{\mathbb R}}

\def\ls{\leqslant}
\def\gs{\geqslant}

\newcommand{\B}{{\mathscr {B}}}
\newcommand{\abs}[1]{\left\vert#1\right\vert}

\begin{document}

\title{\bf Time regularity of L\'{e}vy-type evolution in Hilbert spaces and of some $\alpha$-stable processes.}
\author{{Witold Bednorz \thanks{Institute of Mathematics, University of Warsaw, ul. Banacha 2, 02-097 Warsaw, Poland} \thanks{ Research partially supported by National Science Centre, Poland grant 2016/21/B/ST1/01489.} \  \& Grzegorz G\l{}owienko $\phantom{l}^{*\dagger}$ \& Anna Talarczyk$\phantom{l}^*$ \thanks{ \hbox{e-mail:}
annatal@mimuw.edu.pl. Research supported in part by  National Science Centre, Poland, grant  2016/23/B/ST1/00492.} }
}

\maketitle

\begin{abstract}
In this paper we consider the existence of weakly c\`adl\`ag versions of a solution to a linear
equation in a Hilbert space $H$, driven by a Levy process taking values in
a Hilbert space $U$. In particular we are interested in diagonal type processes, where  process on coordinates  are  functionals of independent  $\alpha$ stable symmetric process.   We give the if and only if characterization in this case. We apply the same techniques to obtain a sufficient condition for existence of a c\`adl\`ag versions of stable processes described as integrals of deterministic functions with respect to symmetric  $\alpha$-stable random measures with $\alpha\in[1,2)$.

\end{abstract}


\medskip


\section{Introduction}

A L\'{e}vy type evolution equation can be formulated as
\begin{align}
dX_t=AX_tdt+dZ_t,\;\;t\in T=[0,a], \label{basic-eq}\;\;X_0=0,\;\;a>0, 
\end{align}
where $X=(X_t)_{t\in T}$ has values in a separable Hilbert space $H$, $A$ is a generator of a $C_0$ semigroup $(S(t))_{t\ge 0}$ on $H$
and $Z=(Z_t)_{t\in T}$ is a L\'{e}vy-type process. This equation has been considered in several papers, see e.g. \cite{16}, \cite{16-5}, \cite{18-5} and some references therein. We refer to \cite{18} for the general theory of stochastic equations in Hilbert spaces with L\'evy noise.

Note that we have not precised in which space $Z$ should
take its values.  It is far from being trivial since in general $Z$ may have values in a much larger Hilbert space $U$
than $H$ whereas still $X$ is well defined in $H$. The problem is well described in the introduction to \cite{18-5}. 
Note that in general the equation (\ref{basic-eq}) has the solution
\begin{align}
X_t=\int^t_0 S(t-s)dZ_s,\;\;X_0=0 \label{basic-solv}.
\end{align}

In the present paper we consider
only the diagonal case with negative diagonal operator $A$ and diagonal L\'{e}vy-type process $Z$, which is a much simpler question. 
Namely, let  $(e_n)^{\infty}_{n=1}$ be an orthonormal and complete basis in $H$, we assume that for any $n=1,2,\ldots$ vector $e_n$ belongs to the domain of $A$ and $Ae_n=-\gamma_n e_n$ with $\gamma_n>0$. Moreover, assume that $Z_t=\sum_{n=1}^\infty Z_t^{(n)}e_n$, where $Z^{(n)}$ are real-valued independent symmetric L\'evy processes without Gaussian part and with L\'evy measures $\mu_n$, respectively. Note that, in general the sum defining $Z$ may not converge in $H$, but in some larger space $U$.
By
the solution to the diagonal type evolution equation we mean the process
$$
X_t=\sum^{\infty}_{n=1} X^{(n)}_t e_n,\;\;t\in [0,1],
$$
where
\begin{equation}
dX^{(n)}_t=-\gamma_n X^{(n)}_t dt+dZ^{(n)}_t,\;\;X^{(n)}_0=0,\;\;n=1,2,\ldots.
\label{e:Xneq}
\end{equation}
The process $X$ takes values in $H$ if and only if the series $\sum_{n=1}^\infty (X_t^{(n)})^2$ converges in probability (and therefore almost surely, thanks to independence). One can write appropriate conditions in terms of L\'evy measures $\mu_n$ (see Proposition 2.6 in \cite{18-5}).
An important example considered in literature is when $Z^{(n)}=\sigma_n L^{(n)}$, where $L^{(n)}$ are independent standard symmetric $\alpha$-stable L\'evy processes and $\sigma_n\ge 0$. This will be referred  to as the $\alpha$-stable case. In this case, the condition for $X$ to take values in $H$ is 
\begin{equation}
 \sum_{n=1}^\infty \frac {\sigma_n^\alpha}{1+\gamma_n}<\infty.
\label{e:XHvalued}
 \end{equation}

The question we treat in this paper is
the regularity of paths of $X$. Obviously, one may think of the 
existence of a c\`{a}dl\`{a}g version of $X$ in $H$. This case is
described in Liu Zhai \cite{16} and the point is that if such a version exists, then $Z$ takes values in $H$ (in the $\alpha$-stable case it is equivalent to $\sum_{n=1}^\infty\sigma_n^\alpha<\infty$).
However, the intriguing situation is when the condition fails. That
means $Z$ has values beyond $H$, but still we expect some regularity 
of $X$. In this paper  we focus on the existence of a cylindrical c\`{a}dl\`{a}g  modification.

According to the Definition 1.1 in \cite{18-5} an $H$-valued process $X$ is cylindrical c\`adl\`ag if for any $z\in H$ the real valued process
\begin{equation}
Y_t=\langle z,X_t\rangle=\sum_{n=1}^\infty\langle z,e_n\rangle X^{(n)}_t, \qquad t\ge 0.
\label{e:Y_def}
\end{equation}
has a   c\`{a}dl\`{a}g modification.
The consequence of the fact is that for any finite set of vectors $z_1,z_2,\ldots,z_n\in H$ the process
$$
(\langle z_1,X\rangle,\langle z_2,X\rangle,\ldots ,\langle z_n,X\rangle)
$$
has a c\`{a}dl\`{a}g modification, which may indicate a good behavior
of $X$ as a process in high dimensions. There are some partial results towards
the question discussed in the extensive paper \cite{18-5}
(note that there are discussed many forms of regularity). 
However, the results in \cite{18-5} do not completely cover even 
 the basic question of $Z^{(n)}$ that are $\alpha$-stable, where $\alpha \in (1,2)$.  We do propose an approach which in particular covers the question
 formulated as Question 4 in \cite{18-5}. It should be mentioned that the case of $\alpha\in (0,1]$ was completely solved in \cite{16-5}. 
 As it will be proved, our approach works in much general setting of diagonal type evolution equation implying a nice sufficient condition for cylindrical  c\'{a}dl\'{a}g property for all diagonal type evolution equations. Therefore, we partially answer also the Question 3 in \cite{18-5}.

 The process $Y$ of \eqref{e:Y_def} clearly depends on $z$, but as most of the time we will work with fixed $z$,  we do not stress this dependence. Note that even if $X$ does not take values in $H$ can still make sense at least for some $z\in H$ and we may consider the problem of a c\`adl\`ag modification.

 The key idea in the proof is to use the Poissonian representation of L\'evy processes and an application of a result of \cite{B-M} concerning suprema of Bernoulli processes. In this approach it is important that the L\`evy processes are symmetric.
\smallskip

 In the last part of the paper we show the usefulness of our method beyond the
evolution equations. Namely, we give sufficient conditions for existence of  c\`adl\`ag modifications of 
stable processes of the form 
\begin{equation}
 X_t=\int_Ef(t,x)M(dx)\qquad t\in [0,a],
 \label{e:stable}
\end{equation}
where $M$ is a symmetric $\alpha$-stable random measure and 
 $f$ is a deterministic function satisfying appropriate integrability conditions. See Section \ref{section:stable} and Theorem \ref{thm:stable_proc} below. It is worth stressing that our condition also works in the case $\alpha\in(1,2)$, which seems to be a difficult one.
 
 \smallskip
 The paper is organized as follows. In Section 2 we introduce some notation and representations of the process $Y$ given by \eqref{e:Y_def}. In Section 3 we discuss a necessary condition for existence of a c\`adl\`ag modification of the process $Y$. In Section 4 we provide a sufficient condition. Finally, in Section 5 we discuss the problem of c\`adl\`ag modification of stable processes of the form \eqref{e:stable}.
 
 \section{Representation of solution}

For the sake of simplicity we assume that $T=[0,1]$.
As we have explained the solution to the evolution equation has the form (\ref{basic-solv}). Suppose that
$Z^{(n)}=\sigma_n L^{(n)}$, where $\sigma_n\gs 0$ and $L^{(n)}$, $n=1,2,\ldots$ are independent symmetric L\'evy processes without Gaussian component and 
with L\'evy measures $\nu_n$, respectively. 
That is, $L_t^{(n)}$ has characteristic function of the form
\begin{equation*}
 \E e^{i\theta L_t^{(n)}}=\exp\left\{-t\int_\R (1-\cos(\theta y))\nu_n(dy)\right\},
\end{equation*}
where $\nu_n$ is a symmetric Borel measure on $\R$, satisfying $\nu_n(\{0\})=0$ and
\begin{equation*}
 \int_\R (y^2\wedge 1) \nu_n(dy)<\infty.
\end{equation*}
Such processes have c\`adl\`ag modification, and in the sequel  we will always assume that $L^{(n)}$, $n=1,2,\ldots$ are c\`adl\`ag. As described in the introduction we assume that $A$ is a diagonal operator, and for an orthonormal basis $(e_n)_{n}$ of $H$ we have 
 $Ae_n=-\gamma_n e_n$,  with $\gamma_n>0$. Then (\ref{e:Xneq}) reads as
\begin{equation}
X^{(n)}_t=\int^t_0 \exp\left(-\gamma_n(t-s)\right) \sigma_ndL^{(n)}_s.
\label{e:Xn}
\end{equation}
It is well known that the jump times and sizes of $L^{(n)}$ are points of a Poisson random measure, with intensity measure $\ell\otimes \nu_n$, where $\ell$ is the Lebesgue measure on $\R_+$. We denote this random measure by $\pi_n$.  Thus 
$$
L^{(n)}_t=\lim_{\delta\to 0}\int^t_0 \int_{\abs{y}
\ge \delta} y{\pi}_n(ds,dy),
$$ 
where the limit is a.s. Moreover, on a subsequence $\delta_n\searrow 0$ fast enough the convergence is a.s. uniform on bounded intervals (see e.g. Theorem 6.8 in\cite{18}). Note that here we do not need to compensate, since $\nu_n$ are symmetric.
\smallskip

\noindent
 Also, due to symmetry $\pi_n$ can be represented as
\begin{equation*}
 \pi_n =\sum_{i}\delta_{(t_{n,i},(\tilde\varepsilon_{n,i}y_{n,i}))},
\end{equation*}
where $(t_{n,i},y_{n,i})$ are points of a Poisson random measure with intensity
$\ell\otimes \mu_n$ with $\mu_n(B)=2\nu_n(B\cap \R_+)$, which will be denoted here by $\pi_n^+$, and $\tilde\varepsilon_{n,i}$ $i=1,2,\ldots$ are i.i.d.  Rademacher random variables. In this setting
the process $L^{(n)}$ at time $t_{n,i}$  has a jump of absolute value $y_{n,i}$ and sign $\tilde\varepsilon_{n,i}$, i.e.
\begin{equation*}
 \Delta L^{(n)}_{t_{n,i}}=\tilde\varepsilon_{n,i}y_{n,i}.
\end{equation*}
For $n=1,2,\ldots$ the corresponding Poisson random measures $\pi_n^+$ and random signs are independent. 
\smallskip

\noindent
An important example is when $L^{(n)}$ are symmetric $\alpha$-stable processes. 
In this case it is well  known that 
$$
\nu_n(dy)=\frac{C_{\alpha}}{|y|^{\alpha+1}}dy.
$$
Here $C_\alpha>0$ is a constant that  standardizes $L^{(n)}$, so that
\begin{equation*}
 \E e^{i\theta L_t^{(n)}}=e^{-t\abs{\theta}^\alpha}.
\end{equation*}

We fix $z\in H$ and consider existence of a c\`adl\`ag modification of 
\begin{equation}
 Y_t=\langle X_t,z\rangle=\sum_{n=1}^\infty Y_t^{(n)}=\sum_{n=1}^\infty \langle z, e_n\rangle X_t^{(n)}, \qquad t\in[0,1].
 \label{e:Yagain}
\end{equation}
where $X^{(n)}$ are given by \eqref{e:Xn}, and $Y_t^{(n)}=\langle z, e_n\rangle X_t^{(n)}$.


Under a weak assumption the sum $\sum_n Y^{(n)}_t$ 
converges a.s. for all $t\in T=[0,1]$, we explain it below. Each of the variables $Y^{(n)}_t$, $t\in T$ can be represented in terms of the Poisson random measure $\pi_n$ as
 \begin{equation}
  Y^{(n)}_t=\lim_{\delta\to 0+}\sum_{i:y_{n,i}\ge \delta}b_n \va_{n,i}y_{n,i}e^{-(t-t_{n,i})\gamma_n}1_{t_{n,i}\ls t},\label{e:zXt}
 \end{equation}
 where  $b_n=|\sigma_n\langle z, e_n\rangle|$, $\varepsilon_{n,i}=\tilde \varepsilon_{n,i}\sgn(\langle z,e_n\rangle)$ and $t_{n,i}, y_{n,i}$, $\tilde\varepsilon_{n,i}$, $n=1,2,\ldots,$ $i=1,2,\ldots$ are as above. 
 
 We have
  \begin{prop}
For any $t>0$ the sum on the right hand side of \eqref{e:Yagain} converges almost 
surely if and only if
\begin{equation}
\psi(\theta):= \sum_{n=1}^\infty \int_0^t\int_{\R}\left(1-\cos\left(\theta b_nye^{-\gamma_n s}\right)\right) \nu_n(dy)ds<\infty, \qquad \theta\in \R
 \label{e:convergence}
\end{equation}
and the function $\psi$ is continuous at $0$.
\end{prop}
This result follows directly from the fact that $Y^{(n)}$ can be written in the form of integrals with respect to compensated Poisson random measure and their independence.

In particular, if $L^{(n)}$ are standard symmetric $\alpha$-stable L\'evy processes, then
\begin{equation*}
 \int_{\R}\left(1-\cos\left(\theta b_nye^{-\gamma_n s}\right)\right) \nu_n(dy)=
\abs{\theta}^\alpha \left(b_nye^{-\gamma_n s}\right)^\alpha
\end{equation*}
and
the series in \eqref{e:Yagain}
converges almost surely for any $t>0$ if and only if 
\begin{equation*}
 \sum_{n=1}^\infty \frac{b_n^\alpha}{1+\gamma_n}<\infty.
\end{equation*}
 
 \smallskip
 
 It is clear that each of the processes $Y^{(n)}$ is c\`adl\`ag. Thus, using \eqref{e:zXt}
 we can use the following representation of $Y$ 
\begin{equation}
Y_t=\langle z, X_t\rangle =\sum_n Y^{(n)}_t=\sum_{n} \sum_{i} b_n \va_{n,i}y_{n,i}e^{-(t-t_{n,i})\gamma_n}1_{t_{n,i}\ls t},\;\;t\in T,
\label{e:sumXn}
\end{equation}
The sum over $i$ is understood as $\lim_{\delta\to 0}\sum_{i: y_i\ge \delta}...$.
We are ready to discuss the convergence of $\sum_n Y^{(n)}_t$, $t\in T$.

The main idea we follow is that $(Y_t)_{t\in T}$ can be split into two parts
according to whether $b_ny_{n,i}\ge 1$ or $b_ny_{n,i}< 1$. The first part is a finite sum of c\`adl\`ag processes  and in the second the series with respect to $n$,
converges uniformly in $L^1$, thus there is a subsequence on which the convergence is a.s. uniform on $T$, hence the limit is c\`adl\`ag.

\section{Necessary condition}
Recall \eqref{e:sumXn} and \eqref{e:convergence}. The next theorem provides a necessary condition for $Y$ to have a c\`adl\`ag modification. This result follows from Theorem 3.4 of \cite{16-5}, but, as it is short, we will also present its proof, to have a full picture of our problem.
\begin{theorem}
\label{thm:necessary}
 If $Y$ has a c\`adl\`ag modification, then for any $\varepsilon>0$ we have
 \begin{equation}
  \sum_{n=1}^\infty\nu_n\left( \left[\frac \varepsilon b_n, \infty\right) \right)<\infty.
  \label{e:tail_nu_n}
 \end{equation}
\end{theorem}
\begin{ex} (Cf. Corollary 3.5 in \cite{16-5}).
 If $L^{(n)}$ are independent standard symmetric $\alpha$-stable L\'evy processes and the process $Y$ has a c\`adl\`ag modification, then
 \begin{equation}
  \sum_n b_n^{\alpha}<\infty.
  \label{e:nec}
 \end{equation}
\end{ex}
\begin{proofof}{Theorem \ref{thm:necessary}}
We argue by contradiction. Suppose that  \eqref{e:tail_nu_n} does not hold for some $\varepsilon>0$ and that $Y$ has a c\`adl\`ag modification $\tilde Y$. Fix any $n$ and denote:
\begin{equation*}
 Y^{(n,\varepsilon)}_t=\sum_{i:y_i\ge \varepsilon}b_n \va_{n,i}y_{n,i}e^{-(t-t_{n,i})\gamma_n}1_{t_{n,i}\ls t}, \qquad t\ge 0.
\end{equation*}
Then the processes 
\begin{equation}
 \tilde Y_t-Y_t^{(n,\varepsilon)}, \ t\ge 0,\qquad \textrm{and}\qquad Y_t^{(n,\varepsilon)}\, t\ge 0
\label{e:proc}
 \end{equation}
are c\`adl\`ag and they are independent (independence follows from the fact that $\pi_n$ is independently scattered). Moreover, $Y^{(n,\varepsilon)}$ has jumps at jump times of the Poisson process $\pi_n([0,t]\times \{y:\abs{y}\ge \varepsilon\}),\ t\ge 0$.  Therefore, with probability one, the sample paths of the two processes defined in \eqref{e:proc} must have jumps at different times. Hence, with probability one, whenever  $Y^{(n)}$ has a jump of size $\ge \varepsilon$, then $\tilde Y$ has a jump of equal size and sign. Notice also, that
\begin{equation*}
\abs{ \Delta {Y^{(n)}_s}}=b_n\abs{\Delta L^{(n)}_s},
\end{equation*}
Where, for a c\`adl\`ag process $Z$ we denote $\Delta Z_s=Z_s-Z_{s-}$. 
\smallskip

\noindent
We will show that if \eqref{e:tail_nu_n} does not hold then, with probability one, there are infinitely many $n$, such that $L^{(n)}$ has a jump of size $\ge \varepsilon/b_n$. Moreover, all $L^{(n)}$ are independent, hence they jump at different times. Consequently, by the argument above, this implies that $\tilde Y$ must have an infinite number of jumps of size $\ge \varepsilon$ on $[0,1]$, and therefore cannot be c\`adl\`ag. This is a contradiction.
\smallskip

\noindent
Let $\xi^{(n)}$ denote the maximal jump of $L^{(n)}$ on $[0,1]$; $\xi^{(n)}=\sup_{s\le 1}\abs{\Delta L_s}$.
Clearly, for $u>0$
$$
\P(\xi^{(n)}< u)=\P({\pi}^{(n)}([0,1]\times \{y:\;\;|y|>u\})=0)=
\exp\left(-\nu_n(\{y:\abs{y}\ge u\})\right).
$$
Hence
 \begin{align}
 \sum_n& \P(b_n \xi^{(n)}\ge \varepsilon)=\sum_n \P(\xi^{(n)}\ge \frac{\varepsilon}{b_n})\notag\\
 &=\sum_n \left(1-\exp(-2\nu_n([\frac \varepsilon{b_n},\infty))\right)
 \label{e:maxjump}\\
 & \gs e^{-1}\sum_n \min\{2\nu_n([\frac \varepsilon{b_n},\infty),1 \}=\infty,\notag
 \end{align}
where the last equality is a consequence of \eqref{e:tail_nu_n}.
As $\xi^{(n)}$ are independent, the 
 Borel Cantelli lemma implies that with probability $1$ there are infinite number of $n$ such that $L^{(n)}$ has a jump of size at least $\varepsilon/b_n$.
\end{proofof}

\section{Sufficient condition}

We now discuss sufficient conditions for existence of c\`adl\`ag modification of $Y$. 
\begin{theorem} \label{thm:sufficient}
 Assume that there exists $\varepsilon>0$ such that \eqref{e:tail_nu_n} is satisied, and additionally that 
 \begin{equation}
  \label{e:sufficient_square}
  \sum_{n=1}^\infty\int_\R b_n^2\int_{b_n\abs{y}\le \varepsilon}\abs{y}^2\nu_n(dy)<\infty.
 \end{equation}
Then $Y$ has a c\`adl\`ag modification.
\end{theorem}
Before we go to the proof of the theorem we make several observations:
\begin{rem}
 The assumptions of Theorem \ref{thm:sufficient} may be also written in the form
 \begin{equation*}
  \sum_{n=1}^\infty \int_\R(\abs{b_ny}^2\wedge 1) \nu_n(dy)<\infty
 \end{equation*}
thus our result is stronger than Theorem 3.8 in \cite{16-5}, where  $\abs{b_ny}$ appeared with power $1$ instead of the square.
 \end{rem}

\begin{ex}
 If $L^{(n)}$ are independent standard symmetric $\alpha$-stable L\'evy processes with $\alpha\in(0,2)$
 then \eqref{e:tail_nu_n} and \eqref{e:sufficient_square} both reduce to 
 \begin{equation}
  \sum_n b_n^{\alpha}<\infty.
  \label{e:necc}
 \end{equation}
 Hence by Theorems \ref{thm:necessary} and \ref{thm:sufficient} \eqref{e:necc} is a necessary and sufficient condition for $Y$ to have a c\`adl\`ag modification. This strengthens the result of \cite{16-5} (Theorem 3.9) which was only proved there for $\alpha<1$.
 \end{ex}
 
 \begin{col}
  Assume \eqref{e:XHvalued}. Then  $X=(X_t)_{t\in T}$, $T= [0,1]$ has cylindrical c\`adl\`ag property if and only if
  \begin{equation}
   \sum_{n=1}^\infty \sigma_n^{\frac{2\alpha}{2-\alpha}}<\infty.
   \label{e:cyl_cadlag_prop}
  \end{equation}
 \end{col}
Recalling the definition of we see
 $b_n$,  \eqref{e:necc} is equivalent to 
 $$
\sum_n |\langle z,e_n\rangle\sigma_n|^{\alpha}<\infty.
$$ For $X$ to have the cylindrical c\`adl\`ag property, \eqref{e:necc} has to be satisfied for any $z\in H$. 
 therefore the Corollary follows by Hahn Banach theorem.
%
%

Note that it is possible that \eqref{e:XHvalued} is satisfied and $\sum_n \sigma_n^{\alpha}=\infty$ but \eqref{e:cyl_cadlag_prop} is satisfied. This means that in this case the process $X$  is not $H$-c\`adl\`ag but it is cylindrically c\`adl\`ag, and for which the process $Z$ of \eqref{basic-eq} does not have values in $H$.

\begin{proofof}{Theorem \ref{thm:sufficient}}
As in the proof of Theorem \ref{thm:necessary} let $\xi^{(n)}$ denote the maximal size of a jump of $L^{(n)}$ on $[0,1]$. Then, by \eqref{e:maxjump} and an elementary estimate $1-e^{-x}\le x$ we have that 
\begin{equation*}
 \sum_n \P(b_n \xi^{(n)}\ge \varepsilon)<\infty.
\end{equation*}
 Borel Cantelli lemma and the fact that each $L^{(n)}$ is c\`adl\`ag imply that there are only a finite number of $y_{n,i}$ such that 
$b_n y_{n,i}\ge \varepsilon$. 
\smallskip

\noindent
Instead of $Y$ it is therefore enough to consider the process
\begin{equation}
 Y^{(\varepsilon)}_t:=\sum_{n=1}^\infty Y_t^{(n,\varepsilon)},\qquad t\ge 0,
 \label{e:Yepsilon}
\end{equation}
where 
\begin{equation}
 Y_t^{(n,\varepsilon)}=\lim_{\delta\to 0}\sum_{i: \delta\le y_{n,i}<\varepsilon}
 b_n\varepsilon_{n,i}y_{n,i}e^{-\gamma_n(t-t_{n,i})}1_{t\ge t_{n,i}},
 \label{e:Xepsilon}
\end{equation}
since  the difference between $Y$ and $Y^{(n,\varepsilon)}$ is a finite sum of c\`adl\`ag processes.
Note that 
\begin{equation*}
  Y_t^{(n,\varepsilon)}=\sigma_n \langle z, e_n\rangle\int_0^t e^{-\gamma_n(t-s)}dL_s^{(n,\varepsilon)},
\end{equation*}
where $L_t^{(n,\varepsilon)}=L_t-\sum_{s\le t: b_n\abs{\Delta_sL}\ge \varepsilon}\Delta L_s$. Each of the processes $Y^{(n,\varepsilon)}$ is c\`adl\`ag.
\smallskip

\noindent
Moreover observe that thanks to \eqref{e:sufficient_square} the process
\begin{equation*}
 L^{(\varepsilon)}=\sum_{n=1}^\infty \sigma_n\langle z, e_n\rangle  L^{(n,\varepsilon)}
\end{equation*}
is well defined and the sum converges in $L^2$ in the supremum norm on $[0,1]$, since $L^{(n,\varepsilon)}$ are independent martingales and 
\begin{equation*}
 \sum_{n=1}^\infty \E (L^{(n,\varepsilon)}_1)^2 =\sum_{n=1}^\infty b_n^2\int_{b_n\abs{y}\le \varepsilon}\abs{y}^2\nu_n(dy)<\infty,
\end{equation*}
by assumption \eqref{e:sufficient_square}. Therefore $L^{(\varepsilon)}$ is c\`adl\`ag.
\smallskip

\noindent
The problem thus reduces to showing that 
\begin{equation}
 L_t^{(\varepsilon)}-Y_t^{(\varepsilon)}=\sum_n\left( L_t^{(n,\varepsilon)}-Y_t^{(n,\varepsilon)}\right),\qquad t\ge 0
 \label{e:Zepsilon}
\end{equation}
has a c\`adl\`ag modification.
\smallskip
 
\noindent
We will show that with probability one the series in \eqref{e:Zepsilon} converges a.s. in the supremum norm. 
The property implies the existence  of a c\`adl\`ag modification of the limit.
Since we could not find the right reference we give a short proof below for the sake of completeness .
  \begin{lem}\label{zbiega}
  Suppose that real processes $(\eta^{(n)}_t)_{t\in T}$, $T=[0,a]$ are independent and  c\`adl\`ag. Moreover, suppose that for any $\va>0$
  \begin{align}
  \lim_{N\ra \infty}\sup_{n\gs m\gs N}\P\left( \left\|\sum^n_{k=m} \eta^{(k)}\right\|_{\infty}>\va\right)=0. \label{cond-sup}
  \end{align}
  Then, for any $t\in [0,1]$ the process $\eta_t=\sum^{\infty}_{n=1}\eta^{(n)}_t$ has a c\`adl\`ag modification. More precisely, $\sum^{\infty}_{n=1}\eta^{(n)}$
  converges a.s. in the Skorohod $J_1$ topology to some $\bar{\eta}$ which the c\`adl\`ag  modification of $Y$. Moreover, the series $\sum^{\infty}_{n=1}\eta^{(n)}$ also converges uniformly.
  \end{lem}
\begin{rem}
Note that the space $D([0,1])$ equipped with the supremum norm is not separable, so we cannot follow the usual approach for separable Banach spaces. In fact we even do not know whether $\omega\ra\sum^m_{n=1}\eta^{(n)}_t$ is a random variable wit values in $D([0,1])$ equipped with the
$\sigma$-field generated by the supremum norm.
\end{rem}
\begin{proofof} {Lemma \ref{zbiega}}
For $x,y\in D([0,1])$ let 
$$
d(x,y)=\inf_{\lambda \in \Lambda}\max\left( \sup_{0\ls s<t\ls 1} \frac{\log[\lambda(t)-\lambda(s)]}{t-s},\|x-y\circ \lambda\|_{\infty} \right),
$$
where $\Lambda$ is the set of nondecreasing continuous functions from $[0,1]$ onto itself. It is known that $d$ is a metric on $D([0,1])$
inducing the Skorohod $J_1$ topology and such that the space $D([0,1])$ with this metric is a Polish space (see \cite{BIL}). Clearly,
$d(x,y)\ls \|x-y\|_{\infty}$, hence
$$
\sup_{n\gs m\gs N}\P\left( d(\sum^n_{k=1} \eta^{(k)},\sum^{m-1}_{k=1}\eta^{(k)} )>\va\right)\ls \sup_{n\gs m\gs N}
\P\left( \left\|\sum^n_{k=m} \eta^{(k)}\right\|_{\infty}>\va \right).
$$
The space $(D([0,1]),d)$ is complete and that is why the series $\sum^{\infty}_{n=1}\eta^{(n)}$ converges in probability in this space.
 By Theorem 1 \cite{KAL} it also converges almost surely in the metric $d$ to some $\bar{\eta}$ which is c\`adl\`ag.
  Moreover, a simple consequence of (\ref{cond-sup}) is that
 $\|\eta^{(n)}\|_{\infty}$ converges in probability to $0$ as $n\ra \infty$. Therefore, by Theorem 2 of \cite{KAL}, the series $\eta=\sum^{\infty}_{n=1}\eta^{(n)}$
 also converges a.s. in the uniform norm.  Therefore, for any fixed $t\in [0,1]$ variables $\eta_t=\bar{\eta}_t$ a.s. 
It completes the proof.
\end{proofof}
The processes $\eta^{(n)}=L^{n,\varepsilon}-Y^{(n,\varepsilon)}$ areindependent for $n=1,2,\ldots$ and  c\`adl\`ag,
therefore  it suffices to prove that the supremum norms converge in $L^2$.

We will prove the following lemma
\begin{lem}\label{lem:cauchy} 
There exists a universal positive constant $C$ such that 
for any $k\le m$ we have
\begin{equation}\E \sup_{t\in[0,1]}\abs{\sum_{n=k}^m \left(L_t^{(n,\varepsilon)}-Y^{(n,\varepsilon)}_t\right)}^2\le C_1 E \abs{\sum_{n=k}^m \left(L_1^{(n,\varepsilon)}-Y^{(n,\varepsilon)}_1\right)}^2\le C_2 \sum_{n=k}^m b_n^2\int_{b_n\abs{y}\le \varepsilon}\abs{y}^2\nu_n(dy).
\label{e:sup}
\end{equation}
\end{lem}
By assumption \eqref{e:sufficient_square} this implies the Cauchy condition for   the series in \eqref{e:Zepsilon}. The proof of the theorem will be complete provided that we show Lemma \ref{lem:cauchy}, which we do presently.\end{proofof}

\begin{proofof}{Lemma \ref{lem:cauchy}}

Denote 
\begin{equation*}
 a_{n,i}(t)=b_ny_{n,i}(1-e^{\gamma_n(t-t_{n,i})})_+.
\end{equation*}
 Then for fixed $k\le m$
 \begin{equation}
  \sum_{n=k}^m \left(L_t^{(n,\varepsilon)}-X^{(n,\varepsilon)}\right)
  =\lim_{\delta\to 0+} A^{(\delta)},
\label{e:sumA}
  \end{equation}
 where for $\delta<\varepsilon$ 
 \begin{equation}
A^{(\delta)}_t =   \sum_{n=k}^m\sum_{i:\delta\le b_n y_{n,i}<\varepsilon}
  \varepsilon_{n,i}a_{n,i}(t).
  \label{e:Adelta}
 \end{equation}
In \eqref{e:sumA} the limit is in $L^2$ for any fixed $t\in [0,1]$ moreover, it is a.s. uniform on $[0,1]$ on a subsequence $\delta_n\searrow 0$ fast enough.
\smallskip

\noindent
We will estimate the expectation of the supremum norm of $A^{(\delta)}$ on $[0,1]$ using  a result of \cite{B-M}. Observe that the double sum in \eqref{e:Adelta} is a.s. finite and the random processes $a_{n,i}$ are nondecreasing, $a_{n,i}\le b_ny_{n,i}$, moreover $(\varepsilon_{n,i})_{n,i}$ are independent of $(a_{n,i})_{n,i}$. The latter processes depend only on $\pi_n^+$, $n=1,2,\ldots.$
Conditioning on $\pi_n^+$ $n=k,\ldots, m$  and using Theorem 1 of \cite{B-M} for any $u> 0$ 
we have
\begin{equation*}
 \P_{\ccE}(\sup_{t\in [0,1]}A^{(\delta )}_t\ge 8 u )\le 53 \P_\ccE(A^{\delta}_1\ge u)
\end{equation*}
 Here $\P_\ccE$ indicates integration with respect to $\varepsilon_{n,i}$ only.
 Taking expectation, using the identity $E\xi^2=2\int_0^\infty uP(\abs{\xi}\ge u)du$ and also symmetry we obtain:
  \begin{equation*}
   \E \sup_{t\in [0,1]}
   \abs{A^{(\delta)}_t}^2 \le C E\abs{A^{(\delta)}_1}^2=
  C \E \sum_{n=k}^m \sum_{i:\delta<b_n\abs{y_{n,i}}\le \varepsilon} b_n^2y_{n,i}^2 a^2_{n,i}(1)
   \le
   \E \sum_{n=k}^m\int_{\delta \le\abs{b_ny}<\varepsilon}b_n^2y^2\nu_n(dy)
  \end{equation*}
  Letting $\delta\to 0$ we obtain \eqref{e:sup}.
\end{proofof}

\section{C\`adl\`ag modification of processes expressed as integrals with respect symmetric stable random measures.}
\label{section:stable}

A large class of stable stochastic processes studied in literature are of the form
\begin{equation}
 X_t=\int_Ef(t,x)M(dx)\qquad t\in [0,a]
 \label{e:stable_proc}
\end{equation}
where $a>0$, $M$ is an $\alpha$-stable random measure defined on some measurable space $(E,\B)$ and $f:[0,a]\times E\mapsto \R$ is a measurable function on the product space, satisfying appropriate integrability conditions. See e.g. \cite{ST} for a systematic treatment of stable integrals and stable processes. In this section we discuss a sufficient condition for the process of the form \eqref{e:stable_proc} to have a c\`adl\`ag modification (and hence for local boundedness of the process). Necessary and sufficient conditions for sample boundedness of processes of the form \eqref{e:stable_proc} in the case  $\alpha<1$ are known. The case $\alpha>1$ seems to be more difficult (see Chapter 10 of \cite{ST}). Some more recent results on c\`adl\`ag property of stable integrals of the form \eqref{e:stable_proc} can be found in
 \cite{DD} and \cite{BR}.

It turns out that our methods used in the previous section can be applied also in this setting in case where $M$ is a symmetric $\alpha$-stable random measure.

We assume that $0<\alpha<2$ and let $m$ be a $\sigma$-finite measure on a measurable space $(E,\B)$. Let $M$ denote a symmetric $\alpha$-stable random measure on $E$ with control measure $m$. That is, if we denote by ${\cal E}_0:=\{A\in \B: m(A)<\infty\}$ then $(M(A))_{A\in {\cal E}_0}$ is a family of real valued random variables such that:
\begin{itemize}
 \item [(i)] For any $A_1,A_2, \ldots \in {\cal E}_0$ such that $A_i\cap A_j=\emptyset$ for $i\neq j$ the random variables $M(A_1),M(A_2),\ldots$ are independent. Moreover, if  we also have that  $m(\bigcup_{n=1}^\infty A_n)<\infty$, then
 \begin{equation*}
  M(\bigcup_{n=1}^\infty A_n)=\sum_{n=1}^\infty M(A_n), \qquad a.s.
 \end{equation*}
 \item[(ii)] If $A\in {\cal E}_0$, then $M(A)$ is a symmetric $\alpha$-stable random variable with scale parameter $(m(A))^{\frac 1\alpha}$, that is
 \begin{equation*}
  \E e^{i\theta M(A)}=\exp\{-m(A)\abs{\theta}^\alpha\}, \qquad \theta \in \R.
 \end{equation*}
\end{itemize}
If $g:E\mapsto \R$ is a measurable function such that
\begin{equation*}
 \int_E\abs{g(x)}^\alpha m(dx)<\infty
\end{equation*}
then one can define $\int_E f(x)M(dx)$. This is done in the usual way, by approximating $f$ by simple functions and passing to the limit. It turns out that for integrals defined in this way we have
\begin{equation*}
 \E \exp\{i \int_Eg(x)M(dx)\}=\exp\{-\int_E\abs{g(x)}^\alpha m(dx)\}.
\end{equation*}
Therefore, if $a>0$, $f:[0,a]\times E\mapsto \R$ is measurable with respect to the  $\sigma$-fields $\B([0,a])\otimes \B /\B(\R)$ and such that for any $t>0$ we have
\begin{equation*}
 \int_{E}\abs{f(t,x)}^\alpha m(dx)<\infty,
\end{equation*}
then the process \eqref{e:stable_proc} is well defined.
\smallskip

\noindent
Recall also, that $M(A)$ and $\int_E f(t,x)M(dx)$ may be constructed using a Poisson random measure. Assume that $\pi$ is a Poisson random measure on $\R\times E$ with intensity measure 
\begin{equation}
 \frac {c_\alpha}{\abs{z}^{1+\alpha}}dz m(dx),
\label{e:falpha}
 \end{equation}
where $c_\alpha>0 $ is chosen such that 
\begin{equation*}
\int_{\R}(1-\cos z)\frac{c_\alpha}{\abs z^{1+\alpha}}dz=1.
\end{equation*}
Then, for $A\in \mathcal{E}_0$ 
\begin{equation*}
 M(A)=\lim_{\delta\to 0}\int_{\{z:\abs{z}>\delta\}\times A}z\pi(dzdx),
\end{equation*}
where the limit is in probability, and a.s. if taken over a sequence $\delta_n\searrow 0$.
\smallskip

\noindent
If $\delta_n\searrow 0$ and $E_n\in \B$ are such that  $m(E_n)<\infty$, $E_n\subset E_{n+1}$ for all $n$ and $\bigcup_n E_n=E$, then for fixed $t$,
the stable integral with respect to the stable random measure constructed above may be represented as
\begin{equation}
 \int_E f(t,x)M(dx)=\lim_{n\to \infty}\int_{\{z:\abs{z}>\delta_n\}\times E_n}zf(t,x)\pi(dzdx), \qquad a.s. \label{e:4.2}
\end{equation}
A simple, but key observation in our context is that since the L\'evy measure $\frac{c_\alpha}{\abs{z}^{1+\alpha}}dz$ is symmetric, the Poisson random measure $\pi$  may be written as
\begin{equation}
 \pi=\sum_{i}\delta_{(\varepsilon_i y_i,x_i)},
\label{e:pi_again}
 \end{equation}
where $\pi^+=\sum_i\delta_{(y_i,x_i)}$ is a Poisson random measure with intensity measure $\frac {2c_\alpha}{y^{\alpha+1}}1_{y>0}dym(dx)$
and $\varepsilon_1,\varepsilon_2,\ldots$ are i.i.d Rademacher random variables independent of $\tilde \pi$.
\smallskip

\noindent
We have the following theorem.
\begin{theorem}\label{thm:stable_proc}
 Assume that $(X_t)_{t\in[0,a]}$ is of the form \eqref{e:falpha} and  $f=f_1-f_2$, where the functions $f_1,f_2:[0,a]\times E\mapsto \R_+$, $i=1,2$  are $\B([0,a])\otimes \B /\B(\R)$ measurable and such that there exists a set $N\in \B$,  $m(N)=0$ such that for any $x\in E\backslash N$ the functions $t\mapsto f_{i}(t,x)$ are c\`adl\`ag and nondecreasing, $i=1,2$. Moreover, assume that 
 \begin{equation}
  \int_E\abs {f_i(a,x)}^\alpha m(dx)<\infty, \qquad i=1,2.
 \label{e:ass_fi}
 \end{equation}
Then the process  $(X_t)_{t\in[0,a]}$ defined by \eqref{e:stable_proc} has a c\`adl\`ag modification.
\end{theorem}
\begin{rem}
 Assumptions of Theorem \ref{thm:stable_proc}  essentially mean that for any $x\in E\backslash N$ the function $t\mapsto f(t,x)$ is c\`adl\`ag and has finite variation on $[0,a]$. Moreover, this variation as a function of $x$ is in $L^\alpha(E, m)$. 
\end{rem}
\begin{proofof}{Theorem \ref{thm:stable_proc}}
 Let $\pi$ be a Poisson  random measure of the form \eqref{e:pi_again} and let $\delta_n$ and $E_n$ be as in \eqref{e:4.2}. 
 Note that $\pi$ restricted to the set 
 $ \{\abs{z}:z>\delta_n\}\times E_n$  is such that the number of 
  points $(\varepsilon_iy_i,x_i)$ in this set is
  Poisson with parameter $\int_{\{y:y>\delta_n\}}\frac{2 c_\alpha}{y^{1+\alpha}}dy\, m(E_n)$ 
 and then all random variables  $\varepsilon_i,y_i,x_i$ $i=1,2,..$ are independent, $\varepsilon_i$ are Rademacher random variables, $y_i$ have law with the density proportional to   $ 1_{(\delta_n,\infty)}(y)\frac 1{y^{1+\alpha}} $ and $x_i$ have the law $\frac 1{m(E_n)}m\big|_{E_n}$.
\smallskip
 
 \noindent
 Let us denote 
 \begin{equation*}
  X_t^{(n)}=\int_{\{\abs{z}>\delta_n\}\times E_n}zf(t,x)\pi(dzdx)=\sum_{i:y_i>\delta_n, x_i\in E_n}\varepsilon_i y_if(t,x_i).
 \end{equation*}
Clearly the process $(X^{(n)}_t)_{t\in [0,a]}$ is c\`adl\`ag since the sum is finite and the function $t\mapsto f(t,x)$ is c\`adl\`ag for any $x\in E\backslash N$. For any $t\in[0,a]$ $X_t^{(n)}$ converges pointwise to $X_t$. Therefore, to prove the theorem it suffices to show that the processes $X^{(n)}$ converge a.s. uniformly on $[0,a]$. Moreover, writing 
\begin{equation*}
 X_t^{(n)} =\int_{\{\abs{z}>\delta_n\}\times E_n}zf_1(t,x)\pi(dzdx)
 -\int_{\{\abs{z}>\delta_n\}\times E_n}zf_2(t,x)\pi(dzdx)
\end{equation*}
it suffices to show that each of the two processes on the right hand side converges a.s. uniformly on $[0,a]$.
\smallskip

\noindent
Hence, without loss of generality in what follows we will assume  that $f=f_1$, i.e. $f$ is nonnegative, $t\mapsto f(t,x)$ is c\`adl\`ag and  nondecreasing for any $x\in E\backslash N$ and $f$ satisfies \eqref{e:ass_fi}.
\smallskip

\noindent
Let us denote
\begin{equation*}
 B_a=\{(z,y)\in \R\times E: \abs{zf(a,x)}\le 1 \}.
\end{equation*}
Thanks to the assumption \eqref{e:ass_fi} it is immediate to see that 
\begin{equation*}
 \int_{B_a^c}\frac{c_\alpha}{\abs{z}^{1+\alpha}}dzm(dx)<\infty,
\end{equation*}
Hence $\pi$ has a finite number of points in $B_a^c$. It is therefore enough to consider only the part of $X^{(n)}$ which is an integral over the set $A_n:=(\{\abs{z}>\delta_n\}\times E_n)\cap B_a$.
\smallskip

\noindent
Denote 
\begin{equation*}
 Y_t^{(n)}:=\int_{A_n}zf(t,x)\pi(dzdx).
\end{equation*}
We will show that 
\begin{equation}
 \lim_{m,n\to \infty}\E\sup_{t\in [0,a]}\abs{Y_t^{(n)}-Y_t^{(m)}}^2=0
\label{e:YnCauchy}
 \end{equation}
This will imply that $Y^{(n)}$ converge in probability in the supremum norm, but since $Y^{(k)}-Y^{(k-1)}$, $k=1,2,\ldots$ are independent 
we can once again use Lemma \ref{zbiega}, which implies that $Y^{(n)}$ converge a.s.  in the supremum norm, thus the limit is c\`adl\`ag.
\smallskip

\noindent
Hence to complete the proof of the theorem it suffices to show \eqref{e:YnCauchy}. This is similar to the proof of Lemma \ref{lem:cauchy}.
Suppose that $n\ge m$, then 
\begin{equation*}
 Y_t^{(n)}-Y_t^{(m)}=\sum_{i:(y_i,x_i)\in (A_n\backslash A_m)}\varepsilon_i y_if(t, x_i)
\end{equation*}
Integrating out first with respect to $\varepsilon_i$ and applying Theorem 1 of \cite{B-M} we have that
\begin{equation*}
 \E\sup_{t\in [0,a]}\abs{Y_t^{(n)}-Y_t^{(m)}}^2\le C
 \E\sum_{i: (y_i,x_i)\in A_n\backslash A_m}y_i^2 f^2(a,x)=\int_{A_n\backslash A_m}y^2f^2(a,x)\frac {2c_\alpha}{y^{\alpha+1}}dy\, m(dx)\to 0.
\end{equation*}
The last convergence follows from the fact that 
\begin{equation*}
 \lim_{n\to \infty}\int_{A_n}y^2f^2(a,x)\frac {2c_\alpha}{y^{\alpha+1}}dy\, m(dx)=\int_{B_a}y^2f^2(a,x)\frac {2c_\alpha}{y^{\alpha+1}}dy\, m(dx)
\end{equation*}
which is finite by assumption \eqref{e:ass_fi}.
 \end{proofof}

\end{document}